\documentclass[11pt,reqno]{amsart}
\usepackage[T1]{fontenc}
\usepackage[utf8]{inputenc}
\usepackage[english]{babel}
\usepackage{lmodern}
\usepackage[a4paper,margin=27mm,headheight=14pt]{geometry}
\usepackage{mathtools,amssymb,amsthm,booktabs,array}
\usepackage{enumitem,needspace}
\usepackage[protrusion=true,expansion=false]{microtype}
\usepackage{xcolor,xurl}
\usepackage[colorlinks=true,linkcolor=blue!40!black,citecolor=green!30!black,urlcolor=blue!40!black]{hyperref}
\hypersetup{pdftitle={Unbounded Nilpotency Gaps on Rationally Elliptic Manifolds},pdfauthor={Paul-Eugène Parent},pdfsubject={Homology nilpotency, rational category, elliptic Sullivan models, and smooth realization}}
\setlist{itemsep=3pt,topsep=4pt}
\AtEndDocument{\raggedbottom}
\allowdisplaybreaks[2]
\newtheorem{theorem}{Theorem}[section]
\newtheorem{proposition}[theorem]{Proposition}
\newtheorem{lemma}[theorem]{Lemma}
\newtheorem{corollary}[theorem]{Corollary}
\theoremstyle{definition}
\newtheorem{definition}[theorem]{Definition}
\theoremstyle{remark}
\newtheorem{remark}[theorem]{Remark}
\newtheorem{example}[theorem]{Example}
\newcommand{\Q}{\mathbb Q}
\newcommand{\N}{\mathbb N}
\newcommand{\Z}{\mathbb Z}
\newcommand{\m}{\mathfrak m}
\newcommand{\Hnil}{\operatorname{Hnil}}
\newcommand{\nilh}{\operatorname{nil}_{h}}
\newcommand{\catz}{\operatorname{cat}_{0}}
\newcommand{\Catz}{\operatorname{Cat}_{0}}
\newcommand{\cat}{\operatorname{cat}}
\newcommand{\Cat}{\operatorname{Cat}}

\newcommand{\wt}{\operatorname{wt}}
\newcommand{\len}{\operatorname{len}}

\newcommand{\norm}[1]{\lVert #1\rVert_1}
\newcommand{\Regr}{\mathcal R}

\newcommand{\fundingmark}{\textsuperscript{*}}
\numberwithin{equation}{section}
\title[Nilpotency gaps on elliptic manifolds]{Unbounded Nilpotency Gaps\protect\\ on Rationally Elliptic Manifolds}
\author[Paul-Eug\`ene Parent]{Paul-Eug\`ene Parent\fundingmark}
\thanks{\makebox[0pt][r]{\fundingmark\enspace}This research was supported by the Natural Sciences and Engineering Research Council of Canada (NSERC), grant no.~149521.}
\address{Department of Mathematics and Statistics, University of Ottawa, Ottawa, Ontario, Canada}
\email{pparent@uottawa.ca}
\date{October 14, 2026}
\subjclass[2020]{55P62, 55M30, 57R65}
\keywords{Homology nilpotency, rational Lusternik--Schnirelmann category, elliptic manifold, minimal Sullivan model, acyclic ideal, rational surgery}
\begin{document}
\begin{abstract}
We construct simply connected, closed smooth manifolds that are rationally
elliptic and for which homology nilpotency exceeds rational category by an
arbitrarily large amount. For every $q\geq5$, a one-parameter family
of even regradings of an explicit minimal Sullivan model with $q+4$
generators admits smooth realizations $V_{q,K}$ satisfying
\[
 e_0(V_{q,K})=\catz(V_{q,K})=\Catz(V_{q,K})=\nilh(V_{q,K})=3q-2,
 \qquad 4q-6\leq\Hnil(V_{q,K})\leq4q-4.
\]
The lower bound follows from two forced primitives and an explicit
orientation functional. The even regrading separates the degrees involved in this obstruction
and makes middle cohomology vanish in the dimensions relevant to
rational surgery. Every resulting model is realized in its formal
dimension, with vanishing rational Pontryagin classes. Rational
connectivity can be made arbitrarily large. Under the spherical stabilization
statement recalled in Section~\ref{sec:stabilization}, a single cone attachment
can increase homology nilpotency arbitrarily while preserving rational category.
\end{abstract}
\maketitle
{\small\tableofcontents}

\section{Introduction}\label{sec:intro}

The rational homotopy type of a simply connected closed manifold has a
minimal Sullivan model, unique up to isomorphism. Other commutative
differential graded models may have smaller nilpotency length. We show
that requiring a model to be a \emph{quotient of the minimal model} can
increase the least possible nilpotency length by an arbitrarily large
amount, even for rationally elliptic manifolds.

For a minimal Sullivan algebra $M$, let $\nilh(M)$ be the least nilpotency
length of any commutative differential graded algebra quasi-isomorphic to
$M$. Let $\Hnil(M)$ be the least such length among quotients $M/J$ for
which $J$ is an acyclic differential ideal. Equivalently,
\begin{equation}\label{eq:intro-Hnil}
 \Hnil(M)=\min\{n:\ \m_M^{n+1}\subset J,\quad
 J\text{ an acyclic differential ideal of }M\}.
\end{equation}
All ideals in this formula are graded by ordinary cohomological
degree. No auxiliary homogeneity is imposed. For a space $X$, both
invariants are evaluated on its minimal rational model. Cornea's
nil-length theorem identifies $\nilh(X)$ with rational strong category
$\Catz(X)$ \cite[Definition~2, Theorem~1.1 and Corollary~2, p.~837]{Cornea1994}.
Thus the distinction is between the
existence of a short model and the existence of a short quotient of the
minimal model.

\subsection{The geometric result}

For each $q$, our main examples form a one-parameter family of even
regradings of the original model $M_q$ introduced in Section~\ref{alg:family}.
We denote the regraded models by
\begin{equation}\label{eq:intro-model}
 \widetilde M_{q,K}
 =\bigl(\Lambda(x,z,a,b_0,\ldots,b_q),D\bigr),
 \qquad q\geq5,
\end{equation}
with two even generators $x,z$ and $q+2$ odd generators
$a,b_0,\ldots,b_q$. Their differentials have the uniform form
\begin{equation}\label{eq:intro-differential}
 Dx=Da=0,\qquad Dz=ax,\qquad
 Db_i=x^{q-i}z^i+iab_{i-1},\qquad b_{-1}=0.
\end{equation}
The parameter $K$ indexes changes of cohomological degrees, not
realizations of a single fixed graded model. For sufficiently large $K$, the regrading separates a finite
range of auxiliary weights. All generator degrees and sufficient bounds on $K$ are given
explicitly in \eqref{alg:constants}--\eqref{alg:newdegrees}.

\begin{theorem}\label{thm:main}
For every $q\geq5$ and every integer $K\geq3W_q+1$, where
$W_q=3q^2+3q-3$, there exists an oriented, simply connected,
closed smooth manifold $V_{q,K}$ whose minimal rational model is
$\widetilde M_{q,K}$. It is rationally elliptic and satisfies
\begin{equation}\label{eq:main-invariants}
 \begin{gathered}
 e_0(V_{q,K})=\catz(V_{q,K})=\Catz(V_{q,K})
 =\nilh(V_{q,K})=3q-2,\\
 4q-6\leq\Hnil(V_{q,K})\leq4q-4.
 \end{gathered}
\end{equation}
Its dimension is the formal dimension $D_{q,K}$ in
\eqref{alg:newtop}, and the realization can be chosen with
$p_i(TV_{q,K})=0$ in rational cohomology for all $i\geq1$.
Moreover,
\[
 \dim_{\Q}\bigl(\pi_*(V_{q,K})\otimes\Q\bigr)=q+4,
\]
and its rational connectivity tends to infinity with $K$.
\end{theorem}

In particular, for each $r\geq1$ the choice $q=r+4$ gives
\begin{equation}\label{eq:prescribed-gap}
 \nilh(V_{r+4,K})=3r+10,\qquad
 r\leq\Hnil(V_{r+4,K})-\nilh(V_{r+4,K})\leq r+2.
\end{equation}
These bounds place $\Hnil$ in an interval of width two, which
suffices to establish unbounded additive gaps among rationally elliptic
closed manifolds. The parameter $q$ controls the length of the generator chain in a
single explicit model; it does not count factors in a product.

The smooth realization is part of the main conclusion. It preserves
the entire prescribed rational homotopy type and its formal dimension.
The manifolds arise from rational surgery applied directly to the
elliptic models. Their construction uses neither a doubling operation
nor the addition of further minimal generators.

There is also a consequence for cone attachments. Under the spherical
stabilization statement of
\cite[Theorem~1.2, p.~8; Theorem~15.4 and Corollary~15.6,
pp.~130--131]{Parent2026}, one can choose
a degreewise finite wedge $S$ of simply connected rational spheres such
that, with $Y=(V_{q,K})_{\Q}\vee S$,
\[
 \Hnil(Y)=\catz(Y)=3q-2.
\]
The homotopy cofibre of $S\hookrightarrow Y$ is $(V_{q,K})_{\Q}$.
Consequently, attaching the single cone $CS$ increases $\Hnil$ by an
amount between $q-4$ and $q-2$, while $\catz$ remains unchanged.
Section~\ref{sec:stabilization} proves this assertion with the required
finiteness qualifications.

\subsection{Literature and the position of the result}

The comparison starts with Lusternik--Schnirelmann category and strong
category. We use their normalized conventions, in which a sphere has
category one. Ganea established the basic comparison between these
invariants \cite{Ganea1967}. Cornea identified strong category with cone
length \cite{Cornea1995} and showed that, in the rational setting, cone
length equals the least nilpotency length of a commutative differential
graded model \cite{CorneaCone1994,Cornea1994}. These results give a geometric meaning
to the flexible invariant $\nilh$. The general theory and its relations
with cone decompositions are developed in \cite{CLOT2003}.
The classical quantitative comparison is
\begin{equation}\label{eq:classical-category-gap}
 \catz(X)\leq\Catz(X)\leq\catz(X)+1;
\end{equation}
see~\cite[Proposition~3.15, p.~82]{CLOT2003}.
The first inequality can be strict: Dupont's counterexample to the
Lemaire--Sigrist conjecture has rational category three and rational
cone length four \cite{Dupont1999}; see also the discussion in
\cite[p.~744]{CFL1998}, which cites Dupont's 1996 preprint.
Thus the classical category defect is bounded
by one, whereas the quotient defect studied here is unbounded.

Sullivan's minimal models provide the algebraic framework
\cite{Sullivan1977}. F\'elix and Halperin described rational category
through the word-length quotients of a minimal model and their homotopy
retractions \cite[Theorem~4.7, p.~16]{FH1982}. Toomer's invariant originates in the Moore
spectral sequence \cite{Toomer1974}; its minimal-model formulation
detects how long an essential cocycle can be. The standard inequalities
are
\begin{equation}\label{eq:comparison-intro}
 e_0(X)\leq\catz(X)\leq\Catz(X)=\nilh(X)\leq\Hnil(X).
\end{equation}
In particular, a homotopy retraction of a word-length projection and a
quasi-isomorphism onto a short quotient express different requirements.
For the quotient condition, all required primitives must lie in a
single ideal, closed under multiplication by every element of the
minimal model.

Poincar\'e duality makes the present setting especially informative.
F\'elix, Halperin, and Lemaire proved that $e_0=\catz$ for simply connected rational
Poincar\'e duality spaces \cite[Theorem~3, p.~750]{FHL1998}. Cornea, F\'elix,
and Lemaire proved that rational category and rational cone length
coincide for simply connected Poincar\'e complexes \cite[Theorem~1.1, p.~744]{CFL1998}.
Together these theorems identify the first four terms of
\eqref{eq:comparison-intro} in our family. The unbounded gap therefore occurs in the last inequality, while
category and strong category agree.

Ellipticity adds a finiteness condition on rational homotopy as well
as on cohomology. Halperin's finiteness theorem establishes
Poincar\'e duality for elliptic minimal models
\cite[Theorem~3, p.~188]{Halperin1977}; see also \cite[Section~32; Proposition~38.3, p.~512]{FHT2001}.
In the present construction both finiteness and duality can be checked
directly: two regular eliminations produce a finite-dimensional quotient
with an explicit perfect cochain pairing. The same calculation therefore establishes ellipticity and identifies
the orientation used in the nilpotency obstruction. Halperin's theorem
on hyperbolicity of the Poincar\'e pairing also applies in even formal
dimension, since the homotopy Euler characteristic is negative
\cite[Theorem~4, p.~191]{Halperin1977}.
Proposition~\ref{geo:pairing} supplies a separate, direct proof for this
family. For the chosen regradings, a stronger conclusion holds: middle
cohomology vanishes in the surgery dimensions.
Remark~\ref{geo:classical-route} explains how these statements differ.

Homology nilpotency was introduced by Carrasquel--Vera as an invariant
of an ideal in a commutative differential graded algebra
\cite[Definition~20]{CV2015}, in connection with upper bounds for
sectional category. Here it is applied to the augmentation ideal of
the minimal Sullivan model. The use of the minimal model is essential: replacing it by a
nonminimal Sullivan algebra can change the quotient problem. The invariant is also distinct from the classical
homotopical nilpotency of $H$-spaces studied by Berstein and Ganea
\cite{BersteinGanea1961}. Throughout this paper, $\nilh$ means
homotopical \emph{nil-length} in the short-model sense above.

The first strict separation used as a starting point in this project
is the example $\nilh=3<4=\Hnil$ in
\cite[Theorem~12.4 and Corollary~12.5]{Parent2026}; the quotient establishing its upper bound
is treated in Appendix~B of that memoir. Appendix~\ref{app:first-separator}
below recalls the example and its primitive obstruction. The present
family addresses the following geometric question: can
the defect grow without bound while the minimal model has finitely
many generators and is itself the rational model of a closed smooth
manifold? Theorem~\ref{thm:main} answers this question affirmatively.
Its proof uses the explicit chain $b_0,\ldots,b_q$, not a tensor
amplification of the earlier separator.

Auxiliary weights and changes of degree have classical precedents.
Body--Douglas study integer weight decompositions in
\cite[pp.~335--336]{BodyDouglas1978} and explicitly treat positive
weights in \cite[Definition, p.~22, and Theorem~1, p.~23]{BodyDouglas1980}.
Compatible weight decompositions are treated systematically by
Cirici--Saleh \cite[Definition~2.1]{CiriciSaleh2024}.
Touz\'e gives a related regrading formalism, including its signs,
in \cite[Definition~6.1 and Lemma~6.3, pp.~2587--2588]{Touze2014}. We use an additive even shift of the
cohomological degree, so all graded-commutativity and differential
signs are preserved. What matters for the present lower bound is a
specific arithmetic feature, proved here: every degree in the finite
obstruction range contains at most one weight. An arbitrary graded
ideal in the new model must therefore obey the same primitive
constraints as a weight-homogeneous ideal in the original model.

It is useful to distinguish two kinds of realization. Lambrechts and
Stanley prove that a simply connected cdga whose cohomology satisfies Poincar\'e
duality has a quasi-isomorphic differential graded representative
whose underlying algebra satisfies Poincar\'e duality
\cite[Theorem~1.1]{LambrechtsStanley2008}. This is an algebraic
realization theorem; it does not assert that such a representative
is an optimal short quotient of a fixed minimal model. Smooth
realization of a prescribed rational homotopy type is a further
step, governed by the rational surgery theorems of Barge and
Sullivan \cite{Barge1976,Sullivan1977}. In our family the same
degree separation used for the nilpotency obstruction implies
vanishing middle cohomology whenever the formal dimension is
divisible by four. All rational surgery conditions can then be met
with zero Pontryagin classes.

Cornea's stabilization theorem provides another point of comparison:
adjoining a suspension can reduce cone length to category
\cite[Theorem~1.1, p.~97]{CorneaCone1994}; its rational spherical-wedge form
is recalled in \cite[Theorem~4.1, p.~747]{CFL1998}.
The stabilization theorem of~\cite{Parent2026} has a stronger
conclusion: it supplies a short quotient of the fixed minimal model,
and hence controls $\Hnil$. This extra conclusion does not follow
from the classical comparison of category and cone length.
A single cone attachment can increase category by at most one
\cite[Theorem~1.32, p.~16]{CLOT2003}. The theorem on Poincar\'e
complexes in \cite[Theorem~1.2, p.~744]{CFL1998} describes the
particularly rigid behavior of a top-cell attachment. Our application
cones off the whole spherical summand supplied by~\cite{Parent2026}.
Under this stabilization statement, a single cofibration can increase
$\Hnil$ arbitrarily while leaving category unchanged. The
spherical summand may be infinite, and the cone need not be a
single cell.

A complementary dimension--connectivity criterion gives sufficient
conditions for $\Hnil=\nilh$ \cite{ParentDC2026}.
Section~\ref{geo:manifolds-section} compares its explicit threshold
with the dimensions of our manifolds. The family remains outside
that sufficient equality range even as its rational connectivity grows.

\subsection{The mechanism of the proof}

There are two algebraic tasks and one geometric task. First, a finite
quotient supplies an upper bound for $\Hnil$ and a Poincar\'e pairing.
An explicit closed functional on the minimal algebra then computes
$e_0$, hence all the category invariants, from the word length of an
orientation cocycle. Second, if an acyclic ideal contains sufficiently long
words, acyclicity forces a primitive of $az^{s-1}$. Multiplicative
closure then forces a second cycle. Its forced primitive, multiplied
by a suitable monomial, is detected by the orientation functional,
contradicting the acyclicity of the ideal.

Primitives must first be obtained inside the acyclic ideal: extracting
weight components need not preserve it. Degree separation then determines
their weights and formulas. Thus regrading extends the argument to every
graded ideal of $\widetilde M_{q,K}$, as required by the definition of $\Hnil$.

Geometrically, degree separation confines middle cohomology to at most one
weight. In dimensions divisible by four, the orientation weight is not
divisible by two, forcing middle cohomology to vanish by Poincar\'e duality.
Rational surgery then realizes the model. Part~I develops the algebraic
construction and gap estimates; Part~II treats smooth realization and
cone attachments.

\section{Geometric and algebraic preliminaries}\label{sec:preliminaries}

\subsection{Rational spaces, ellipticity, and category}

All spaces under consideration are simply connected and of finite
rational type in each degree. A closed manifold is compact and without
boundary. We write $X_{\Q}$ for the rationalization of $X$ and evaluate
rational invariants of a manifold on its rationalization without
changing the notation. In particular,
$\catz(X)=\cat(X_{\Q})$ and $\Catz(X)=\Cat(X_{\Q})$.
Category is normalized: $\cat(*)=0$, and $\cat(S^n)=1$ for $n\geq1$.
The ordinary category of a space is the least $n$ for which it admits
a cover by $n+1$ open sets contractible within the space. Strong
category is the least $n$ for which some space of the same homotopy
type admits a cover by $n+1$ contractible open sets. Equivalently,
it is cone length \cite{Cornea1995}.

A space is \emph{rationally elliptic} if both $H^*(X;\Q)$ and
$\pi_*(X)\otimes\Q$ are finite-dimensional. In a minimal Sullivan
model $(\Lambda V,d)$, this means that $V$ and $H(\Lambda V,d)$
are finite-dimensional. A rationally elliptic manifold need not
have finitely many nonzero integral homotopy groups. Similarly,
$X$ is \emph{rationally $(p-1)$-connected} if $\pi_j(X)\otimes\Q=0$
for $j<p$; this does not require the corresponding torsion groups
to vanish.

A graded algebra $H$ has \emph{Poincar\'e duality of formal dimension
$N$} if $H^N\cong\Q$, $H^{>N}=0$, and multiplication followed by
an orientation $H^N\to\Q$ gives nondegenerate pairings
\[
 H^j\otimes H^{N-j}\longrightarrow\Q
 \qquad(0\leq j\leq N).
\]
Elliptic minimal algebras have this property
\cite{Halperin1977}. For any simply connected rational Poincar\'e
duality space, the classical comparison theorems give
\begin{equation}\label{eq:PD-equalities}
 e_0(X)=\catz(X)=\Catz(X)=\nilh(X)
 \quad\text{\cite{FHL1998,CFL1998,Cornea1994}}.
\end{equation}
No equality with $\Hnil$ is part of these theorems.

\subsection{Minimal models and the two nilpotency invariants}

We work over $\Q$. A cdga is a connected, augmented commutative
differential graded algebra in nonnegative cohomological degrees;
the differential has degree $+1$. For homogeneous elements,
\[
 uv=(-1)^{|u||v|}vu,
 \qquad d(uv)=d(u)v+(-1)^{|u|}u\,d(v).
\]
The notation $\Lambda V$ denotes the free graded-commutative
algebra, polynomial on even generators and exterior on odd
generators. A simply connected minimal Sullivan model is
$M=(\Lambda V,d)$ with $V^{\leq1}=0$, a Sullivan ordering, and
$dV\subset\Lambda^{\geq2}V$. Existence, uniqueness, comparison,
and spatial realization are used in their standard forms
\cite{Sullivan1977,FHT2001}. If the differential is decomposed by
word-length increase, our convention is
$d_i:V\to\Lambda^{i+1}V$; thus $d_0=0$ in a minimal model and
$d_1$ is the quadratic part.

Set $\m_M=M^+=\Lambda^{\geq1}V$. Its power
$\m_M^k=\Lambda^{\geq k}V$ consists of sums of words of length
at least $k$. Word length is distinct from cohomological degree.
For a connected augmented cdga $B$, its \emph{nilpotency length}
is the least $n\in\N=\{0,1,\ldots\}$ with $(B^+)^{n+1}=0$.

\begin{definition}\label{def:nilpotencies}
For a simply connected minimal Sullivan algebra $M$, define
\begin{align*}
 \nilh(M)&=\inf\{n:\ M\simeq B,\ (B^+)^{n+1}=0\},\\
 \Hnil(M)&=\inf\{n:\ \m_M^{n+1}\subset J\subset M,
   \ J\text{ an acyclic differential ideal}\}.
\end{align*}
Here $\simeq$ means a zigzag of cdga quasi-isomorphisms, every
ideal is cohomologically graded, and the infimum of the empty
set is $+\infty$.
\end{definition}

An ideal $J$ is acyclic when $H(J)=0$ as a complex. It is not
enough for every cycle in $J$ to be a boundary in $M$:
each such cycle must have a primitive in $J$. This distinction drives the
lower-bound argument.

\Needspace*{8\baselineskip}
\begin{proposition}\label{prop:quotient-criterion}
For a minimal algebra $M$, the inequality $\Hnil(M)\leq n$
holds if and only if there is a surjective quasi-isomorphism
$M\twoheadrightarrow C$ with $(C^+)^{n+1}=0$.
Consequently $\nilh(M)\leq\Hnil(M)$, and both invariants
are well defined on rational homotopy types.
\end{proposition}
\nopagebreak[4]
\begin{proof}
The exact sequence $0\to J\to M\to M/J\to0$ shows that the
projection is a quasi-isomorphism exactly when $J$ is acyclic.
An acyclic ideal is proper and lies in $\m_M$, so
$(M/J)^+=\m_M/J$. Its $(n+1)$st power vanishes exactly when
$\m_M^{n+1}\subset J$. Conversely, the kernel of a surjective
quasi-isomorphism has these properties. Such a quotient is
one of the models allowed in the definition of $\nilh$.
Finally, an isomorphism of minimal algebras preserves powers
of the augmentation ideal and acyclic differential ideals.
Uniqueness of the minimal model then gives the claimed invariance.
\end{proof}

The restriction to a minimal model cannot be omitted. If
$M\to B$ is a quasi-isomorphism with $B^1=0$, it becomes
surjective after adjoining, for each element of a homogeneous basis
of $B^+$, a contractible pair $u,du$ mapped to that element and its
differential.
The enlarged, nonminimal Sullivan algebra therefore has a
short quotient whenever $B$ is short. This construction does
not provide a quotient of $M$ itself.

\subsection{The Toomer invariant and orientation functionals}

We use the word-length description
\begin{equation}\label{eq:Toomer}
 e_0(M)=\sup\{k:\text{some }z\in\m_M^k\cap\ker d
                       \text{ has }[z]\ne0\}.
\end{equation}
Equivalently, $e_0$ is the least $n$ for which
$H(M)\to H(M/\m_M^{n+1})$ is injective
\cite[Remark~9.3(1), p.~29]{FH1982}; see also~\cite{FHT2001}.
The following elementary observation
will turn an explicit top-degree functional into an exact
calculation of category.

\begin{lemma}\label{lem:orientation-bound}
Suppose $H(M)$ has Poincar\'e duality of dimension $N$, with
orientation represented by a cocycle $\omega\in\m_M^n$.
If a linear functional $\varepsilon:M^N\to\Q$ satisfies
\[
 \varepsilon d=0,\qquad \varepsilon(\omega)=1,
 \qquad \varepsilon(\m_M^{n+1}\cap M^N)=0,
\]
then $e_0(M)=n$.
\end{lemma}
\begin{proof}
The orientation gives $e_0(M)\geq n$. If a cycle
$z\in\m_M^{n+1}$ represented a nonzero class, duality would
give a closed $v$ with $[zv]$ a nonzero multiple of $[\omega]$.
The functional would then be nonzero on $zv$, contradicting
$zv\in\m_M^{n+1}$.
\end{proof}

There is a second consequence of duality that clarifies
the ideal obstruction. If $J$ is acyclic, the map
$M\to M/J$ preserves the orientation class. A normalized
closed orientation functional can therefore be chosen to
vanish on $J$. Since $J$ is an ideal, this functional also vanishes on every
product $jv$ with $j\in J$. A forced primitive whose
product has a nonzero orientation value contradicts this
requirement. The lower-bound proof makes the functional and both forced
primitives explicit.

\subsection{Weights and even changes of degree}

A weight grading is a decomposition $M=\bigoplus_{\chi\in\N^2}M_\chi$
such that $D(M_\chi)\subset M_\chi$ and
$M_\chi M_\psi\subset M_{\chi+\psi}$. It is auxiliary:
the weights are not rational homotopy degrees. We write
$\norm\chi=\chi_1+\chi_2$. An additive map
$\sigma:\Z^2\to2\Z$ defines a new cohomological grading by
\begin{equation}\label{eq:general-regrading}
 |u|'=|u|+\sigma(\chi)\qquad(u\in M_\chi).
\end{equation}
Throughout this paper an \emph{admissible even regrading} is such a
change for which every new generator degree is at least two. This
condition retains the simply connected minimal-model setting.
The regrading gives another cdga with the same multiplication and
differential.
The shift is even, so parities and all signs are unchanged.
Since $D$ preserves weights, it still has degree $+1$.

Weight-homogeneous maps can be regraded, and taking cohomology one weight
at a time shows that quasi-isomorphisms are preserved. This does not imply
invariance of $\Hnil$, since cohomologically graded ideals may change.
We will control these ideals by degree separation and compute category
from the regraded orientation and functional.

\part{Elliptic models and nilpotency gaps}\label{part:algebra}

\section{The elliptic family and a finite Poincar\'e duality model}
\label{alg:family}

The construction starts with a chain of odd generators attached to two even
ones. The integer $q$ controls both the length of this chain and the gap. In this section the degrees are deliberately small. They make the
calculations transparent; a change of degrees in Section~\ref{alg:regrading}
will ensure that the obstruction applies to every graded differential ideal.

For $q\geq5$, set
\begin{equation}\label{alg:model}
 M_q=\bigl(\Lambda(x,z,a,b_0,\ldots,b_q),D\bigr),\qquad
 |x|=2,\quad |z|=4,\quad |a|=3,\quad |b_i|=2q-1+2i,
\end{equation}
and define
\begin{equation}\label{alg:differential}
 Dx=Da=0,\qquad Dz=ax,\qquad
 Db_i=x^{q-i}z^i+iab_{i-1},\qquad b_{-1}=0.
\end{equation}
Here $\Lambda$ denotes the free graded-commutative algebra: $x,z$ are even
and $a,b_0,\ldots,b_q$ are odd. We order the odd generators as written. The identity $D^2=0$ follows from
\[
 D(x^{q-i}z^i)=iax^{q-i+1}z^{i-1},\qquad
 D(iab_{i-1})=-iax^{q-i+1}z^{i-1}.
\]
The order $x,a,z,b_0,\ldots,b_q$ is a Sullivan order, and every differential
is decomposable. Thus $M_q$ is a simply connected minimal Sullivan algebra
with two even and $q+2$ odd generators.

We use the auxiliary multigrading
\begin{equation}\label{alg:weights}
 \wt(x)=(1,0),\qquad \wt(a)=(0,1),\qquad
 \wt(z)=(1,1),\qquad \wt(b_i)=(q,i).
\end{equation}
Weights add under multiplication and are preserved by $D$. Every generator
has positive weight sum. In particular, each weight component is finite
dimensional. For a monomial $m=x^u z^t a^e b_S$, where $e\in\{0,1\}$,
$S\subset\{0,\ldots,q\}$, and $k=|S|$, write $\wt(m)=(w,g)$. Then
\begin{equation}\label{alg:count}
 \begin{aligned}
 w&=u+t+qk,& g&=t+e+\sum_{i\in S}i,\\
 |m|&=2w+2g+e-k,& \len(m)&=w-(q-1)k+e.
 \end{aligned}
\end{equation}
We shall repeatedly use the consequence
\begin{equation}\label{alg:degreebound}
 0\leq |m|\leq3\norm{\wt(m)}.
\end{equation}
Indeed, $e\leq g$, so $2w+2g+e-k\leq2w+3g$.

\subsection{Eliminating two regular pairs}

The finite model needed below is obtained directly from $M_q$. The following
lemma permits the elimination of a regular pair even when the remaining
generators do not generate a differential subalgebra.

\begin{lemma}\label{alg:regularpair}
Suppose that the underlying graded algebra of a cdga $A$ is
$A_0\otimes\Lambda u$, where $u$ is odd, $Du=f\in A_0$, and $f$ is a
non-zero-divisor in $A_0$. Then the differential ideal $(u,f)$ is acyclic.
No stability of $A_0$ under the differential is required.
\end{lemma}
\begin{proof}
Write $Dv=D_0v+u\partial v$ for $v\in A_0$. Since $D^2u=0$,
the ideal $(u,f)=fA_0\oplus uA_0$ is differential. Let $fv+uv'$ be a
cycle. The component of its differential in $A_0$ is $f(D_0v+v')$.
Regularity of $f$ therefore gives $v'=-D_0v$. Consequently
\[
 fv+uv'=fv-uD_0v=D(uv),
\]
and the primitive $uv$ belongs to the same ideal.
\end{proof}

\Needspace*{11\baselineskip}
\begin{proposition}\label{alg:quotient}
There is a weight-preserving surjective quasi-isomorphism
\begin{equation}\label{alg:finitequotient}
 M_q\longrightarrow Q_q=
 \frac{\Lambda(x,z,a,b_1,\ldots,b_{q-1})}
 {(x^q,z^q+qab_{q-1})}.
\end{equation}
Its kernel is $(b_0,x^q,b_q,z^q+qab_{q-1})$. The vector-space dimension
of $Q_q$ is $q^2 2^q$, and its augmentation ideal has nilpotency length
exactly $4q-4$.
\end{proposition}
\nopagebreak[4]
\begin{proof}
First eliminate $(b_q,Db_q)$. The element
$Db_q=z^q+qab_{q-1}$ is monic as a polynomial in the even variable $z$,
and hence is regular in the algebra generated by the remaining generators.
Lemma~\ref{alg:regularpair} applies. The resulting graded algebra is free
over $\Q[x]$: polynomial division in $z$ supplies a basis with
$z$-exponent less than $q$. The element $x^q$ is therefore still regular,
and the lemma applies a second time to $(b_0,Db_0)=(b_0,x^q)$.
The composite projection is a quasi-isomorphism. All the relations are
homogeneous for the weights~\eqref{alg:weights}.

A basis of $Q_q$ is
\begin{equation}\label{alg:basis}
 x^u z^t a^e b_S,\qquad
 0\leq u,t\leq q-1,\quad e\in\{0,1\},\quad
 S\subset\{1,\ldots,q-1\}.
\end{equation}
There are $q^2 2^q$ such monomials.
To compute nilpotency, consider monomials before reducing powers of $z$.
A nonzero one has $x$-exponent at most $q-1$. If its $z$-exponent is at
most $q-1$, its length is at most
$(q-1)+(q-1)+q=3q-2$. If that exponent is at least $q$, the relation
$z^q=-qab_{q-1}$ prohibits any additional occurrence of $a$ or $b_{q-1}$.
Moreover $z^{2q}=0$. The largest possible length in this case is
\[
 (q-1)+(2q-1)+(q-2)=4q-4.
\]
It is attained, since
\begin{equation}\label{alg:longproduct}
 x^{q-1}z^{2q-1}b_1\cdots b_{q-2}
 =(-1)^{q-1}q\,ax^{q-1}z^{q-1}b_1\cdots b_{q-1}\ne0.
\end{equation}
This proves the asserted nilpotency length.
\end{proof}

The finite number of generators and the finite-dimensional quotient
\eqref{alg:finitequotient} show that $M_q$ is elliptic. Thus the general
elliptic duality theorem applies~\cite{Halperin1977}. In this family one can
also see the duality directly, including its formal dimension. Set
\begin{equation}\label{alg:top}
 T_q=ax^{q-1}z^{q-1}b_1\cdots b_{q-1},\qquad
 d_q=|T_q|=3q^2+2q-2.
\end{equation}
The element $T_q$ spans $Q_q^{d_q}$, and $Q_q^{d_q-1}=0$.
Indeed, every other basis monomial in~\eqref{alg:basis} is obtained by
removing factors of degree at least two from $T_q$.

Let $\tau_q$ extract the coefficient of $T_q$. Multiplication followed by
$\tau_q$ pairs every basis monomial with its complementary monomial.
A product whose $z$-exponents have sum at least $q$ reduces to a monomial
with $z$-exponent at most $q-2$, and cannot contribute to $T_q$.
Thus the matrix of this pairing is a signed permutation matrix.
The functional $\tau_q$ annihilates boundaries, because there are no
cochains in degree $d_q-1$. Closedness of $\tau_q$ makes boundaries orthogonal to cycles;
nondegeneracy and rank--nullity give $(\ker D)^\perp=\operatorname{im}D$.
The induced pairing on $\ker D/\operatorname{im}D$ is therefore perfect,
proving Poincar\'e duality on $H(Q_q)$ in formal dimension $d_q$.

\subsection{An orientation of controlled word length}

Write
\[
 R_j=xb_j-zb_{j-1}\quad(1\leq j\leq q-1),\qquad
 \omega_q=az^{q-1}\prod_{j=1}^{q-1}R_j.
\]
To check that $\omega_q$ is a cocycle, write $D=\delta+a\theta$ on
expressions without $a$, where
\[
 \delta b_i=x^{q-i}z^i,\quad\delta x=\delta z=0,
 \qquad \theta z=x,\quad\theta b_i=ib_{i-1},\quad\theta x=0.
\]
Then $\delta R_j=0$, and multiplication by $a$ kills the remaining
$a\theta$ terms. The image of $\omega_q$ in $Q_q$ is $T_q$: every other
nonzero term in the product contains $b_0$. In particular, $[\omega_q]\ne0$.

For $0\leq i<j\leq q$, put
\begin{equation}\label{alg:mij}
 m_{ij}=ax^{2q-1-i-j}z^{i+j-1}b_{\widehat{i,j}},
 \qquad b_{\widehat{i,j}}=
 \prod_{\substack{0\leq\ell\leq q\\\ell\notin\{i,j\}}}b_\ell,
\end{equation}
with the product taken in increasing order. The expansion of the orientation
is particularly simple:
\begin{equation}\label{alg:orientation}
 \omega_q=az^{q-1}\prod_{j=1}^{q-1}(xb_j-zb_{j-1})
 =\sum_{k=0}^{q-1}(-1)^k m_{kq}.
\end{equation}
In a nonzero term, the choices $-zb_{j-1}$ must precede all choices
$xb_j$; any reversal repeats an odd generator. This gives the displayed
expansion. Every term has word length $3q-2$, original degree $d_q$, and
weight
\begin{equation}\label{alg:gamma}
 \gamma_q=(A_q,B_q),\qquad
 A_q=q^2+q-2,\qquad B_q=\frac{q(q+1)}2.
\end{equation}

\section{A closed trace and the category of the family}
\label{alg:categorysection}

The coefficient trace on $Q_q$ proves duality, but it does not directly
determine the largest $k$ for which the orientation class has a representative
in $\m^k$. For this purpose we construct a second trace, supported entirely on words
of length $3q-2$. Its uniqueness will later turn the same calculation into
an obstruction to short quotients.

The coefficients below solve the closedness recurrence~\eqref{alg:recurrence};
the proof also checks the remaining sources and normalization.
Define $f_q:M_q\to\Q$ to be zero outside degree $d_q$, zero on all
monomials other than the $m_{ij}$, and to satisfy
\begin{equation}\label{alg:trace}
 f_q(m_{ij})=
 \frac{(-1)^{q+i+j-2}(j-i)\binom qi\binom qj}
 {(2q-1)\binom{2q-2}{i+j-1}}.
\end{equation}
All nonzero values therefore occur in weight $\gamma_q$.

\begin{lemma}\label{alg:closedtrace}
The functional $f_q$ satisfies
\[
 f_qD=0,\qquad f_q(\omega_q)=1,\qquad f_q(\m^{3q-1})=0.
\]
\end{lemma}
\begin{proof}
The last assertion follows from the word lengths of the $m_{ij}$.
To prove closedness, weight preservation allows us to consider only
monomials of weight $\gamma_q$ and degree $d_q-1$.
Formula~\eqref{alg:count} gives $k-e=q-1$ for such a monomial. When $e=0$, the possibilities are exactly
\begin{equation}\label{alg:tracefirstsources}
 c_{ij}=x^{2q-2-i-j}z^{i+j}b_{\widehat{i,j}},
 \qquad 0\leq i<j\leq q,\quad i+j\leq2q-2.
\end{equation}
When $e=1$, they are exactly
\begin{equation}\label{alg:tracesecondsources}
 c'_k=ax^{q-1-k}z^{k-1}b_{\widehat k},
 \qquad 1\leq k\leq q-1,
\end{equation}
where $b_{\widehat k}$ omits $b_k$ from $b_0\cdots b_q$.
These lists follow by solving the two weight equations for the exponents
of $x$ and $z$; their displayed ranges are precisely the conditions that
those exponents be nonnegative.

Write $f_{ij}=f_q(m_{ij})$, with $f_{jj}=f_{i,q+1}=0$.
In $Dc_{ij}$, with $c_{ij}$ as in~\eqref{alg:tracefirstsources}, the terms
without $a$ vanish under $f_q$, while those containing $a$ give
\begin{equation}\label{alg:recurrence}
 (i+j)f_{ij}+(i+1)f_{i+1,j}+(j+1)f_{i,j+1}=0.
\end{equation}
To verify this identity from~\eqref{alg:trace}, put $t=i+j$ and
$\Delta=j-i$. After factoring out the common binomial coefficients and
signs, the sum of the two shifted terms has numerator
\[
 (q-i)(\Delta-1)+(q-j)(\Delta+1)=(2q-t-1)\Delta.
\]
The denominators are related by
\[
 \frac{t}{\binom{2q-2}{t-1}}
 =\frac{2q-1-t}{\binom{2q-2}{t}}.
\]

In $Dc'_k$, the term removing $b_\ell$ has sign
$(-1)^{\ell-[k<\ell]+1}$, where $[k<\ell]$ is one or zero according as
$k<\ell$ or not. Substituting~\eqref{alg:trace} reduces the required
identity, after division by the nonzero factor
$(-1)^{q+k-1}\binom qk$, to
\begin{equation}\label{alg:secondidentity}
 \sum_{\ell=0}^q
 \frac{(k-\ell)\binom q\ell}
 {(2q-1)\binom{2q-2}{k+\ell-1}}=0,
 \qquad 1\leq k\leq q-1.
\end{equation}
The summand with $\ell=k$ is zero. An elementary beta integral verifies
this finite rational identity. Namely,
\[
 \frac{1}{(2q-1)\binom{2q-2}{j}}
 =\int_0^1 t^j(1-t)^{2q-2-j}\,dt.
\]
The sum in~\eqref{alg:secondidentity} becomes
\[
 \begin{aligned}
 &\int_0^1 t^{k-1}(1-t)^{q-k-1}(k-qt)\,dt\\
 &\hspace{15mm}=
 \int_0^1\frac{d}{dt}\bigl(t^k(1-t)^{q-k}\bigr)\,dt=0.
 \end{aligned}
\]
Both endpoint values vanish, since $1\leq k\leq q-1$.
Since the two lists exhaust the possible sources, this proves $f_qD=0$.

Finally, the orientation expansion~\eqref{alg:orientation} and the same
integral formula give
\[
 \begin{aligned}
 f_q(\omega_q)
 &=\int_0^1t^{q-1}
   \sum_{k=0}^{q-1}(q-k)\binom qk
          t^k(1-t)^{q-1-k}\,dt\\
 &=\int_0^1qt^{q-1}\,dt=1.
 \end{aligned}
\]
Here $(q-k)\binom qk=q\binom{q-1}k$ reduces the sum to the binomial
expansion of $(t+(1-t))^{q-1}$.
\end{proof}

\begin{proposition}\label{alg:category}
For the original model one has
\begin{equation}\label{alg:categoryvalues}
 e_0(M_q)=\catz(M_q)=\Catz(M_q)=\nilh(M_q)=3q-2.
\end{equation}
The same equalities hold after every additive even regrading induced
by the weights~\eqref{alg:weights} for which all new generator degrees
are at least two.
\end{proposition}
\begin{proof}
Apply Lemma~\ref{lem:orientation-bound} to the essential orientation
$\omega_q\in\m^{3q-2}$ and the closed normalized trace of
Lemma~\ref{alg:closedtrace}. It gives $e_0(M_q)=3q-2$.

For a simply connected rational Poincar\'e duality type, the equalities
$e_0=\catz$ and $\catz=\Catz$ are respectively the results of
F\'elix--Halperin--Lemaire~\cite[Theorem~3, p.~750]{FHL1998} and
Cornea--F\'elix--Lemaire~\cite[Theorem~1.1, p.~744]{CFL1998}. Cornea's identification of rational
cone length with homotopical nil-length gives
$\nilh=\Catz$~\cite{Cornea1994}. This proves~\eqref{alg:categoryvalues}.

For the regrading assertion, note that the quotient, its kernel, the
orientation, and the trace are all homogeneous in weight. An additive even shift of
degrees preserves the differential, multiplication, parity, and word
length. Exactness is preserved separately on each weight component.
The quotient remains finite dimensional, and its top basis monomial remains
$T_q$, now in its new degree. Its complementary-monomial pairing still gives
Poincar\'e duality. The same essential orientation and closed trace
therefore give the same value of the Toomer invariant, and the classical
equalities apply again. No invariance of arbitrary short quotients
under a change of degrees is being assumed.
\end{proof}

\section{The obstruction from two forced primitives}
\label{alg:obstruction}

An acyclic ideal defining a short quotient contains all sufficiently
long words and must also contain primitives of its cycles. In this family,
two successive primitive problems have unique solutions, and the trace
detects a monomial multiple of the second primitive.
We first treat weight-homogeneous ideals; Section~\ref{alg:regrading}
removes this restriction by choosing suitable degrees.

\Needspace*{9\baselineskip}
\begin{lemma}\label{alg:uniquetrace}
Let $3q-1\leq s\leq4q-4$. A functional $f$ supported in original degree
$d_q$ and weight $\gamma_q$, satisfying
\[
 fD=0,\qquad f(\m^s)=0,\qquad f(\omega_q)=1,
\]
is necessarily $f_q$.
\end{lemma}
\begin{proof}
In degree $d_q$ and weight $\gamma_q$, formula~\eqref{alg:count}
gives $k-e=q-2$. If $e=1$, solving the weight equations gives exactly the
monomials $m_{ij}$, of length $3q-2$. If $e=0$, every possible monomial
has length $4q-4$ and is killed by $f$. Closedness on the sources
$c_{ij}$ therefore imposes the recurrence~\eqref{alg:recurrence}.
The only pair absent from the source list is $(q-1,q)$. Since $i+j>0$,
descending induction on $i+j$ determines every value $f(m_{ij})$ from
$f(m_{q-1,q})$. Thus the space of such closed functionals has dimension
at most one before normalization. Lemma~\ref{alg:closedtrace} supplies a
nonzero normalized member and proves the assertion.
\end{proof}

\begin{theorem}\label{alg:homogeneousobstruction}
Let $q\geq5$ and $3q-1\leq s\leq4q-6$.
No acyclic differential ideal of $M_q$ homogeneous in both original degree
and weight contains $\m^s$.
\end{theorem}
\begin{proof}
Suppose that $J$ is such an ideal. The orientation cocycle $\omega_q$
cannot belong to $J+D(M_q)$. Indeed, if $\omega_q=j+Dv$, then $Dj=0$, and
acyclicity of $J$ would make $j$ a boundary in $J$, contradicting
$[\omega_q]\ne0$. Choose a linear functional vanishing on $J+D(M_q)$
and taking value one on $\omega_q$. Homogeneity allows us to restrict it
to degree $d_q$ and weight $\gamma_q$ without losing these properties.
Since $\m^s\subset J$, Lemma~\ref{alg:uniquetrace} identifies this
functional with $f_q$. Because $J$ is an ideal, we obtain
\begin{equation}\label{alg:idealannihilator}
 f_q(jv)=0\qquad(j\in J,\ v\in M_q).
\end{equation}

First consider the cycle
\[
 Z=az^{s-1}\in\m^s\subset J.
\]
It has weight $(s-1,s)$ and degree $4s-1$. Acyclicity supplies a primitive
in $J$ of degree $4s-2$ and weight $(s-1,s)$. Formula~\eqref{alg:count} forces
$k=e$. If $k=e=0$, the weight equations would give $x$-exponent $-1$.
If $k=e=1$, with odd generator $b_i$, that exponent is $i-q$; hence
$i=q$ is the only possibility. The entire candidate space is therefore
\[
 \Q\,az^{s-q-1}b_q.
\]
Since $D(az^{s-q-1}b_q)=-Z$, the forced primitive is
\begin{equation}\label{alg:firstprimitive}
 U=-az^{s-q-1}b_q\in J,\qquad DU=Z.
\end{equation}

Next, $z^s\in J$ and the element
\begin{equation}\label{alg:secondcycle}
 C=z^s-sxU=z^s+saxz^{s-q-1}b_q
\end{equation}
is a cycle in $J$. It has weight $(s,s)$ and degree $4s$.
A primitive has degree $4s-1$, so~\eqref{alg:count} now gives $k=e+1$.
If $e=0$, the same weight equations leave only $z^{s-q}b_q$.
If $e=1$, there are two indices $i<j$, and the $x$-exponent is
$i+j-2q+1$. It is nonnegative only for $(i,j)=(q-1,q)$.
Thus the space of possible primitives is spanned by the two monomials
\[
 z^{s-q}b_q,\qquad az^{s-2q}b_{q-1}b_q.
\]
All their displayed exponents are nonnegative because $s\geq3q-1$.
Their differentials are
\begin{align*}
 D(z^{s-q}b_q)
 &=z^s+(s-q)axz^{s-q-1}b_q+qaz^{s-q}b_{q-1},\\
 D(az^{s-2q}b_{q-1}b_q)
 &=-axz^{s-q-1}b_q+az^{s-q}b_{q-1}.
\end{align*}
These two differentials are linearly independent. Comparing their
coefficients with~\eqref{alg:secondcycle} determines the unique primitive
\begin{equation}\label{alg:secondprimitive}
 P_s=z^{s-q}b_q-qaz^{s-2q}b_{q-1}b_q\in J.
\end{equation}

Finally, choose
\begin{equation}\label{alg:multiplier}
 V_s=x^4z^{4q-6-s}b_0b_1\cdots b_{q-4}.
\end{equation}
The upper restriction $s\leq4q-6$ is exactly what makes the exponent
of $z$ nonnegative. Reordering the odd generators gives
\begin{equation}\label{alg:detectorproduct}
 P_sV_s=(-1)^{q-3}x^4z^{3q-6}b_0\cdots b_{q-4}b_q
        -q\,m_{q-3,q-2}.
\end{equation}
The first summand has length $4q-4$ and is annihilated by $f_q$.
The second is detected by~\eqref{alg:trace}, contradicting
\eqref{alg:idealannihilator} because $P_s\in J$:
\begin{equation}\label{alg:detector}
 f_q(P_sV_s)
 =\frac{(-1)^q q^3(q-1)}
 {2(2q-1)(2q-3)(2q-5)}\ne0.\qedhere
\end{equation}
\end{proof}

The product $P_sV_s$ need not be closed: the contradiction uses the
vanishing of $f_q$ on the whole ideal $J$, not merely on its cocycles.
The nonzero value in~\eqref{alg:detector} is independent of $s$.
The parameter determines the two primitive spaces and ensures that
the multiplier~\eqref{alg:multiplier} belongs to the algebra. The argument reaches
$s=4q-6$; it does not settle the next two thresholds.

\begin{example}[The first member of the family]\label{alg:qfive}
For $q=5$, the category is $13$. At the threshold $s=14$, the two
forced primitives are
\[
 U=-az^8b_5,\qquad P_{14}=z^9b_5-5az^4b_4b_5.
\]
The multiplier is simply $V_{14}=x^4b_0b_1$, and the obstruction reads
\[
 f_5(P_{14}V_{14})=-\frac{50}{63}.
\]
Thus an ideal containing all words of length $14$ cannot be both acyclic
and homogeneous in weight. Section~\ref{alg:regrading} extends this
obstruction to every graded ideal of the regraded model. The same two-primitive
mechanism therefore produces larger gaps as $q$ increases.
\end{example}

\section{Separating degrees and arbitrary acyclic ideals}
\label{alg:regrading}

The preceding argument does not suffice for an arbitrary ideal in $M_q$
with its original grading. A primitive could contain terms of different
weights in the same cohomological degree, and a weight component of an
element of an ideal need not belong to that ideal. We therefore choose new
degrees so that each relevant cohomological degree admits at most one weight.
The primitives will first be obtained inside the given ideal; the degree
separation will then identify them.

\subsection{An explicit choice of degrees}

With $A_q,B_q$ as in~\eqref{alg:gamma}, put
\begin{equation}\label{alg:constants}
 \begin{gathered}
 G_q=A_q+B_q,\qquad W_q=2G_q+1=3q^2+3q-3,\\
 R_q=W_q+1,\qquad L_q=2W_qR_q+1,\\
 \lambda_1=L_q+1,\qquad\lambda_2=L_q+R_q,
 \qquad K\geq3W_q+1.
 \end{gathered}
\end{equation}
For the weight-isolation proof, the requirements on $K$ are only
\begin{equation}\label{alg:K-sufficient}
 2K>3W_q,\qquad 2KL_q>d_q+1.
\end{equation}
The bound $K\geq3W_q+1$ in~\eqref{alg:constants} is a convenient
sufficient choice, not an optimized one.
For a weight-homogeneous cochain of original degree $d$ and weight
$\chi=(w,g)$, prescribe the new degree
\begin{equation}\label{alg:regrade}
 d'=d+\sigma(w,g),\qquad
 \sigma(w,g)=2K(\lambda_1w+\lambda_2g).
\end{equation}
Denote the resulting cdga by $\widetilde M_{q,K}=\Regr_\sigma M_q$.
The shift is additive, so multiplication is homogeneous. It is even, so
all signs are unchanged. Since $D$ preserves weights, it still has degree
one. Set
\begin{equation}\label{alg:parameters}
 p=2+2K\lambda_1,\qquad \eta=3+2K\lambda_2.
\end{equation}
The new generator degrees are
\begin{equation}\label{alg:newdegrees}
 |x|'=p,\qquad |a|'=\eta,\qquad |z|'=p+\eta-1,
 \qquad |b_i|'=qp-1+i(\eta-1).
\end{equation}
The differential remains exactly~\eqref{alg:differential}, and the same
Sullivan order proves minimality.

The finite quotient~\eqref{alg:finitequotient} and its acyclic kernel
are homogeneous in weight. They therefore admit the same change of degrees.
More explicitly, cohomology decomposes as
\begin{equation}\label{alg:regradedcohomology}
 H^n(\Regr_\sigma M_q)=
 \bigoplus_{\chi\in\N^2}
 H^{n-\sigma(\chi)}((M_q)_\chi).
\end{equation}
Thus $\widetilde M_{q,K}$ is elliptic. The regraded quotient has the same
basis and the same nilpotency length as $Q_q$. Its top monomial
is still $T_q$, because every generator has positive degree. Its formal
dimension is
\begin{equation}\label{alg:newtop}
 \begin{aligned}
 D_{q,K}&=d_q+2K(\lambda_1A_q+\lambda_2B_q)\\
 &=A_qp+B_q(\eta-1)+2-q.
 \end{aligned}
\end{equation}
In particular, this is the new degree of $\omega_q$.

\begin{lemma}\label{alg:isolation}
Every cohomological degree at most $D_{q,K}+1$ in
$\widetilde M_{q,K}$ contains cochains of at most one weight and at most
one original degree.
\end{lemma}
\begin{proof}
The constants satisfy
\[
 L_q\leq\lambda_i\leq2L_q,\qquad 2K>3W_q,
 \qquad 2KL_q>d_q+1.
\]
If a monomial of original degree $d$ and weight $\chi$ has new degree
at most $D_{q,K}+1$, positivity gives
\[
 2KL_q\norm\chi\leq d_q+1+4KL_qG_q.
\]
Consequently
\[
 \norm\chi\leq2G_q+\frac{d_q+1}{2KL_q}<W_q.
\]
By~\eqref{alg:degreebound}, $0\leq d<3W_q<2K$.
If two such monomials have the same new degree, reduction modulo $2K$
therefore identifies their original degrees. Their weights $\chi,\chi'$
then satisfy
\[
 L_q(\norm\chi-\norm{\chi'})+
 (\chi_1+R_q\chi_2)-(\chi'_1+R_q\chi'_2)=0.
\]
Both expressions in the last two parentheses belong to $[0,W_qR_q]$,
an interval of length less than $L_q$. It follows that the weight sums
are equal and then that the two expressions themselves are equal.
Every coordinate is less than $W_q<R_q$, so uniqueness of base-$R_q$
expansion gives $\chi=\chi'$.
\end{proof}

\subsection{The bounds for all graded differential ideals}

\begin{theorem}\label{alg:main}
For every $q\geq5$ and every integer $K$ satisfying~\eqref{alg:constants},
the minimal Sullivan algebra $\widetilde M_{q,K}$ is simply connected
and elliptic, and
\begin{equation}\label{alg:mainvalues}
 \begin{gathered}
 e_0(\widetilde M_{q,K})=\catz(\widetilde M_{q,K})
 =\Catz(\widetilde M_{q,K})=\nilh(\widetilde M_{q,K})=3q-2,\\
 4q-6\leq\Hnil(\widetilde M_{q,K})\leq4q-4.
 \end{gathered}
\end{equation}
The lower bound is taken over all cohomologically graded acyclic differential
ideals, without any hypothesis of weight homogeneity.
\end{theorem}
\begin{proof}
Ellipticity was established above. Proposition~\ref{alg:category} gives
the equalities in~\eqref{alg:mainvalues}. Regrading the quotient of
Proposition~\ref{alg:quotient} gives the upper bound for $\Hnil$.
It remains to exclude every acyclic graded differential ideal containing
$\m^{4q-6}$.

Set $s=4q-6$, which satisfies $s\geq3q-1$ precisely when $q\geq5$,
and suppose $J$ is such an ideal in $\widetilde M_{q,K}$.
As in the proof of
Theorem~\ref{alg:homogeneousobstruction}, the essential orientation is
not in $J+D(\widetilde M_{q,K})$. Since $J$ is graded in the new degree,
there is a functional $\varepsilon$ supported in degree $D_{q,K}$,
vanishing on $J$ and on boundaries, with $\varepsilon(\omega_q)=1$.
Lemma~\ref{alg:isolation} identifies that degree with original degree
$d_q$ and weight $\gamma_q$. The functional, read in those coordinates,
is $f_q$ by Lemma~\ref{alg:uniquetrace}. Thus it annihilates every product
of an element of $J$ with any cochain.

We next check that the two forced primitive arguments remain valid for
this ideal. The cycles $Z$ and $C$ used in
Theorem~\ref{alg:homogeneousobstruction} have weights $(s-1,s)$ and $(s,s)$,
respectively; both are bounded coordinatewise by $\gamma_q$, because
\[
 A_q-s=q^2-3q+4>0,\qquad
 B_q-s=\frac{(q-3)(q-4)}2\geq0.
\]
Their original degrees and those of their primitives are at most
$4s\leq d_q$ for $q\geq5$. Positivity of the shift therefore puts all
these new degrees in the range of Lemma~\ref{alg:isolation}.

The cycle $Z=az^{s-1}$ lies in $J$. Acyclicity first gives a homogeneous
primitive inside $J$. Its new degree admits only the weight $(s-1,s)$
and original degree $4s-2$, since the explicit primitive
$-az^{s-q-1}b_q$ already occurs there. The exhaustive enumeration in
Theorem~\ref{alg:homogeneousobstruction} thus forces
$U=-az^{s-q-1}b_q\in J$. Closure under multiplication then puts
$C=z^s-sxU$ in $J$. Acyclicity gives a primitive in $J$, and
Lemma~\ref{alg:isolation} forces this primitive to have original degree
$4s-1$ and weight $(s,s)$.
The second enumeration identifies it as $P_s$ in
\eqref{alg:secondprimitive}. Hence $P_sV_s\in J$, whereas
\eqref{alg:detector} yields
\[
 0=\varepsilon(P_sV_s)=f_q(P_sV_s)\ne0.
\]
This contradiction excludes $\m^{4q-6}\subset J$, which gives
$\Hnil(\widetilde M_{q,K})\geq4q-6$ with the stated indexing convention.
\end{proof}

The first few bounds are as follows; all entries refer to the regraded
models.
\begin{center}
\begin{tabular}{@{}ccccc@{}}
\toprule
$q$ & Minimal generators & $\nilh=\catz$ & $\Hnil$ & Gap at least\\
\midrule
5 & 9 & 13 & $[14,16]$ & 1\\
6 & 10 & 16 & $[18,20]$ & 2\\
7 & 11 & 19 & $[22,24]$ & 3\\
8 & 12 & 22 & $[26,28]$ & 4\\
\bottomrule
\end{tabular}
\end{center}

For later use, put $q=r+4$. The theorem gives
\begin{equation}\label{alg:prescribedgap}
 \begin{gathered}
 \nilh(\widetilde M_{r+4,K})=3r+10,\\
 4r+10\leq\Hnil(\widetilde M_{r+4,K})\leq4r+12,\\
 r\leq\Hnil(\widetilde M_{r+4,K})-\nilh(\widetilde M_{r+4,K})\leq r+2.
 \end{gathered}
\end{equation}
These bounds prove that the gap is unbounded in the elliptic class.
Moreover, $x$ is the generator of least degree, namely
$p=2+2K\lambda_1$. Increasing $K$ therefore makes the rational connectivity
arbitrarily large, without changing any of the invariant bounds.

\part{Smooth realization and cone attachments}\label{part:geometry}

\section{The middle dimension and rational surgery}
\label{geo:surgery-section}

The algebraic construction provides both the nilpotency gaps and the
Poincar\'e duality needed for geometric realization. We now ask whether
these rational homotopy types can be realized by closed smooth manifolds.
Rational surgery answers this question once the obstruction in the middle
dimension has been checked. For the present family, this check reduces to
the parity of the orientation weight.

Recall from \eqref{alg:gamma} and \eqref{alg:newtop} that the orientation has
weight and regraded degree
\begin{equation}\label{geo:orientation-data}
 \begin{gathered}
 \gamma_q=(A_q,B_q),\qquad
 A_q=q^2+q-2,\qquad B_q=\frac{q(q+1)}2,\\
 D_{q,K}=d_q+2K(\lambda_1A_q+\lambda_2B_q),
 \end{gathered}
\end{equation}
where $d_q=3q^2+2q-2$.
The weight decomposition passes to cohomology because the differential
preserves weights. By Lemma~\ref{alg:isolation}, the cochains in any
cohomological degree at most $D_{q,K}+1$ have at most one weight.
Here we use this property only in the middle degree.

\begin{lemma}[An orientation-weight obstruction]\label{geo:middle-lemma}
Let $H$ be a finite-dimensional graded Poincar\'e duality algebra over $\Q$
of even formal dimension $D$.  Suppose that $H$ also has a grading by
$\mathbb Z^s$, compatible with the cohomological grading and
multiplication, and that $H^D=H^D_\gamma$ for some weight $\gamma$. If $H^{D/2}$ is supported in at most one weight
and $\gamma\notin2\mathbb Z^s$, then $H^{D/2}=0$.
\end{lemma}

\begin{proof}
Suppose instead that $0\ne u\in H^{D/2}$.  There is a weight $\alpha$
such that all of $H^{D/2}$ lies in $H_\alpha$.
Poincar\'e duality gives $v\in H^{D/2}$ with $uv\ne0$.
This product is a nonzero multiple of the orientation.  Its weight is
both $2\alpha$ and $\gamma$, contrary to the hypothesis.
\end{proof}

The conclusion concerns the entire middle cohomology group. In the surgery
dimensions, it is stronger than the vanishing of the signature of the middle
form. This distinction matters in rational surgery: a zero signature alone
does not ensure that a rational symmetric form satisfies the arithmetic
condition.

\begin{proposition}[Vanishing in the surgery dimensions]\label{geo:middle}
For every $q\geq5$ and every integer $K\geq3W_q+1$,
\begin{equation}\label{geo:middle-vanishing}
 4\mid D_{q,K}
 \quad\Longrightarrow\quad
 H^{D_{q,K}/2}(\widetilde M_{q,K})=0.
\end{equation}
\end{proposition}

\begin{proof}
If $q$ is odd, then $d_q$ is odd and the regrading adds an even integer.
Thus $D_{q,K}$ is odd.

If $q\equiv0\pmod4$, both $A_q$ and $B_q$ are even.  The added term
in \eqref{geo:orientation-data} is divisible by four, whereas
$d_q\equiv2\pmod4$.  Hence $D_{q,K}\equiv2\pmod4$.

The only remaining possibility is $q\equiv2\pmod4$.  In this case
$B_q=(q/2)(q+1)$ is odd, so $\gamma_q\notin2\mathbb Z^2$.
If $D_{q,K}=4m$, the weight-isolation result applies to degree $2m$.
The finite quotient \eqref{alg:finitequotient} gives Poincar\'e duality,
so Lemma~\ref{geo:middle-lemma} yields $H^{2m}(\widetilde M_{q,K})=0$.
\end{proof}

The argument applies to both parities of $K$.
Indeed, the constants in \eqref{alg:constants} have
$W_q$ odd, $R_q$ even, $L_q$ odd, $\lambda_1$ even and
$\lambda_2$ odd.  For $q\equiv2\pmod4$ this gives
\begin{equation}\label{geo:K-parity}
 D_{q,K}\equiv2+2K\pmod4.
\end{equation}
Thus, in this residue class of $q$, dimensions divisible by four occur
precisely for odd $K$. The realization argument covers these parameters
as well.

We recall a sufficient form of the rational surgery theorem. Its conclusion
preserves the complete rational homotopy type, as required by the invariants
studied here.

\begin{theorem}[Barge--Sullivan, a sufficient case]\label{geo:surgery}
Let $X$ be a simply connected rational space of finite type whose
cohomology is a Poincar\'e duality algebra of formal dimension $n\geq5$.
Suppose that either $4\nmid n$, or $n=4m$ and $H^{2m}(X;\Q)=0$.
Then there exists a closed, oriented, simply connected smooth
$n$-manifold $V$ such that
\[
 V_{\Q}\simeq X,
 \qquad p_i(TV)=0\in H^{4i}(V;\Q)\quad(i\geq1).
\]
\end{theorem}

\begin{proof}
In dimensions at least five, Sullivan realizes a prescribed rational
homotopy type and rational Pontryagin classes by a simply connected
manifold \cite[Section~13, especially Theorem~13.2, pp.~320--322]{Sullivan1977}
when Barge's numerical and Witt-group conditions hold
\cite[Theorem~1, p.~471]{Barge1976}.
Prescribe $p_i=0$. Barge uses stable-normal classes
\cite[p.~470]{Barge1976}; vanishing is equivalent for tangent and normal
classes because $p(TV)p(\nu_V)=1$. If $4\nmid n$, the theorem applies
directly.

If $n=4m$, choose a rational fundamental class $\mu$.
The middle form has rank zero, so its signature is zero and its Witt
class lies in the image of $W(\Z)\to W(\Q)$. Equivalently, its image
under $W(\Q)\to W(\Q,\Z)$ vanishes, as required by Barge
\cite[pp.~470--471]{Barge1976}.
All prescribed Pontryagin numbers vanish, as do those of $S^{4m}$,
and the signature equation is
\[
 \bigl\langle L_m(0,\ldots,0),\mu\bigr\rangle
 =0=\operatorname{sign}(X).
\]
Thus all the rational surgery conditions hold.
\end{proof}

The next two results extend smooth realizability to all admissible even
weight regradings, using hyperbolicity in dimensions divisible by four.
Degree separation remains essential to our proof of the lower bound for
arbitrary graded acyclic ideals; realizability alone does not transfer
that bound between regradings.

\begingroup
\setlength{\parskip}{0pt plus 1pt}
\begin{proposition}[Hyperbolicity of the Poincar\'e pairing]
\label{geo:pairing}
Let $Q$ be the regrading of $Q_q$ induced by an admissible even
weight regrading of $M_q$, and let $n$ be its formal dimension.
If $n$ is even, the Poincar\'e pairing on $H^{\mathrm{even}}(Q)$ is
hyperbolic. In particular, if $n=4m$, the middle form on $H^{2m}(Q)$
has Witt class zero over $\Q$.
\end{proposition}

\begin{proof}
Let $\tau$ be the coefficient of the top monomial $T_q$.
The complementary-monomial pairing $B(u,v)=\tau(uv)$ is perfect.
By the normal basis~\eqref{alg:basis}, $T_q$ is still the unique top
monomial. Since all generator degrees are at least two,
every other basis monomial has degree at most $n-2$, so $Q^{n-1}=0$
and $\tau D=0$.

Even formal dimension is equivalent to even $q$. In this case
the restrictions of $B$ to $Q^{\mathrm{even}}$ and $Q^{\mathrm{odd}}$
are nondegenerate; the former is symmetric and the latter alternating.
Put $E=(aQ)^{\mathrm{even}}\subset Q^{\mathrm{even}}$.
The normal basis involves $q$ odd generators and gives
\[
 \dim Q^{\mathrm{even}}=q^2 2^{q-1},\qquad
 \dim E=q^2 2^{q-2}=\tfrac12\dim Q^{\mathrm{even}}.
\]
Since $a^2=0$, the subspace $E$ is isotropic and hence Lagrangian.
Thus the symmetric cochain pairing is hyperbolic.

Set $I=D(Q^{\mathrm{odd}})$. For odd $u,v$ and even $z$,
closedness of the trace gives
\[
 B(Du,Dv)=\tau D(uDv)=0,\qquad B(z,Du)=-B(Dz,u).
\]
Nondegeneracy on $Q^{\mathrm{odd}}$ implies
\[
 I^\perp=\ker(D:Q^{\mathrm{even}}\longrightarrow Q^{\mathrm{odd}}),
 \qquad I^\perp/I=H^{\mathrm{even}}(Q).
\]
Here orthogonal complements are taken in $Q^{\mathrm{even}}$.
To show directly that the induced form on the quotient is hyperbolic,
we construct a Lagrangian. Write $e=\dim E$, $i=\dim I$ and
$r=\dim(E\cap I)$. Since $E^\perp=E$, the pairing between $E$ and $I$
has rank $i-r$, so $\dim(E\cap I^\perp)=e-i+r$.
Consequently
\[
 L=\bigl((E\cap I^\perp)+I\bigr)/I
 \subset I^\perp/I
\]
is isotropic of dimension $e-i$, half the dimension of the quotient.
It is a Lagrangian for the induced Poincar\'e pairing.

If $n=4m$, complementary even degrees outside degree $2m$ give
hyperbolic orthogonal summands. Removing these summands leaves the
middle form, whose Witt class is therefore zero.
\end{proof}

\begin{corollary}[Realization of all admissible even regradings]
\label{geo:all-regradings}
Let $q\geq5$ and let $\Regr_\sigma M_q$ be any admissible even weight
regrading of $M_q$. Its rational homotopy type is realized, in its
formal dimension, by a closed, oriented, simply connected, rationally
elliptic smooth manifold with all rational Pontryagin classes zero. In particular,
the original model $M_q$ has such a realization in dimension $d_q$.
\end{corollary}
\begin{proof}
The regraded finite quotient supplies Poincar\'e duality, and its top
monomial has degree $n\geq7q-4\geq31$, since each even generator has
degree at least two and each odd generator at least three.
If $4\nmid n$, apply Theorem~\ref{geo:surgery}.
If $n=4m$, Proposition~\ref{geo:pairing} gives zero middle Witt class
and zero signature. Prescribing all Pontryagin classes to be zero
then satisfies Barge's arithmetic condition and signature equation;
the Pontryagin numbers are those of the sphere.
The full Barge--Sullivan theorem therefore realizes the rational type
\cite[Theorem~1, p.~471]{Barge1976}
\cite[Theorem~13.2, pp.~321--322]{Sullivan1977}.
For the original degrees, $d_q$ is odd when $q$ is odd and is
congruent to two modulo four when $q$ is even.
\end{proof}

\endgroup
\Needspace*{14\baselineskip}
\begin{remark}[Classical hyperbolicity and middle vanishing]
\label{geo:classical-route}
Halperin's theorem gives a hyperbolic Poincar\'e pairing on the even
cohomology of a connected elliptic minimal algebra of even formal
dimension and negative homotopy Euler characteristic
\cite[Theorem~4 and the preceding definition, p.~191]{Halperin1977}.
Its hypothesis on the homotopy Euler characteristic holds here because
$\chi_\pi=-q$. Proposition~\ref{geo:pairing} proves the hyperbolicity
directly for this family. Hyperbolicity describes the
bilinear form, while the rational homotopy type remains elliptic.
Halperin's manifold corollary on the same page also yields signature
zero in dimension $4m$: the closed odd generator $a$ represents a
nonzero class, since minimal boundaries are decomposable.
For the selected regradings, Proposition~\ref{geo:middle} gives the
stronger conclusion that the entire middle cohomology group vanishes
in dimensions divisible by four.
We use this conclusion in the main realization theorem below.
The lower bound for arbitrary acyclic ideals still requires the
separated degrees of $\widetilde M_{q,K}$.
\end{remark}

\section{Closed rationally elliptic manifolds}
\label{geo:manifolds-section}

We now apply the surgery criterion to the regradings used in the gap
theorem. The resulting manifolds inherit the algebraic invariants
computed in Part~\ref{part:algebra}.

\begin{theorem}[Smooth realization of the entire family]
\label{geo:realization}
For every integer $q\geq5$ and every integer $K\geq3W_q+1$, there
exists a closed, oriented, simply connected smooth manifold $V_{q,K}$
whose minimal Sullivan model is $\widetilde M_{q,K}$.
It has dimension
\begin{equation}\label{geo:dimension}
 \dim V_{q,K}
 =3q^2+2q-2+
 2K\left(\lambda_1(q^2+q-2)
       +\lambda_2\frac{q(q+1)}2\right),
\end{equation}
and it can be chosen with all rational Pontryagin classes zero.
Moreover, $V_{q,K}$ is rationally elliptic and satisfies
\begin{equation}\label{geo:invariants}
 \begin{gathered}
 e_0(V_{q,K})=\catz(V_{q,K})=\Catz(V_{q,K})
 =\nilh(V_{q,K})=3q-2,\\
 4q-6\leq\Hnil(V_{q,K})\leq4q-4.
 \end{gathered}
\end{equation}
\end{theorem}

\begin{proof}
The Sullivan realization of $\widetilde M_{q,K}$ is a simply connected
rational space $X_{q,K}$ with that minimal model
\cite[Theorem~8.1, p.~301]{Sullivan1977}.
The finite quotient \eqref{alg:finitequotient} and its regrading show that
$H^*(X_{q,K};\Q)$ is finite-dimensional and satisfies Poincar\'e
duality of formal dimension $D_{q,K}$.
In particular, $D_{q,K}\geq d_q\geq83$.

If $4\nmid D_{q,K}$, Theorem~\ref{geo:surgery} applies directly.
If $4\mid D_{q,K}$, Proposition~\ref{geo:middle} verifies its remaining
hypothesis.  We obtain a manifold with
\[
 (V_{q,K})_{\Q}\simeq X_{q,K}
\]
and with zero rational Pontryagin classes.
Uniqueness of minimal Sullivan models identifies its minimal model
with $\widetilde M_{q,K}$.  Formula~\eqref{geo:dimension} is exactly
the formal dimension \eqref{alg:newtop}.

The minimal model has $q+4$ generators: two in even degrees and
$q+2$ in odd degrees.  Consequently
\[
 \dim_{\Q}\bigl(\pi_*(V_{q,K})\otimes\Q\bigr)=q+4<\infty.
\]
Its cohomology is finite-dimensional, as established above.  These
are the two finiteness conditions defining rational ellipticity.

Finally, all the invariants in \eqref{geo:invariants} depend only on
the rational homotopy type.  For $\Hnil$ specifically, an isomorphism
of minimal models identifies their augmentation-ideal powers and
their acyclic differential ideals.  The algebraic estimates of
Theorem~\ref{alg:main} therefore transfer without alteration.
\end{proof}

The manifolds in Theorem~\ref{geo:realization} have the formal dimensions
of the regraded models, given by~\eqref{geo:dimension}. The dimensions
$83,118,159,206$ for $q=5,6,7,8$ belong instead to the realizations of the
original models supplied by Corollary~\ref{geo:all-regradings}.

\Needspace*{10\baselineskip}
The table illustrates the dimensions obtained from our sufficient bounds.
Here $p$ is the first nonzero rational homotopy degree, so the rational
connectivity is $p-1$.
\begin{center}
\begin{tabular}{rrrr}
\toprule
$q$ & $K$ & $p$ & $D_{q,K}$\\
\midrule
5 & 262 & $8\,024\,538$ & $345\,738\,951$\\
6 & 370 & $22\,574\,442$ & $1\,378\,952\,378$\\
6 & 371 & $22\,635\,454$ & $1\,382\,679\,276$\\
\bottomrule
\end{tabular}
\end{center}
The first two rows use $K=3W_q+1$; the third illustrates a dimension
divisible by four. These large dimensions reflect the deliberately
generous separation bounds, not a claim of dimensional optimality.

\Needspace*{10\baselineskip}
\begin{corollary}\label{geo:unbounded}
For every integer $r\geq1$ there is a closed, simply connected,
rationally elliptic smooth manifold $V$ such that
\[
 r\leq\Hnil(V)-\nilh(V)
   =\Hnil(V)-\catz(V)\leq r+2.
\]
These manifolds can be required to have arbitrarily large rational
connectivity and zero rational Pontryagin classes.
\end{corollary}

\begin{proof}
Take $q=r+4$ in Theorem~\ref{geo:realization}.
The two bounds become
\[
 \nilh(V_{r+4,K})=3r+10,
 \qquad
 4r+10\leq\Hnil(V_{r+4,K})\leq4r+12.
\]
It remains to check the connectivity assertion.
The smallest degree of a minimal generator is
\[
 p=|x|=2+2K\lambda_1.
\]
Indeed, $|a|=3+2K\lambda_2>p$, $|z|=p+|a|-1>p$,
and $|b_i|=qp-1+i(|a|-1)>p$.
The correspondence between minimal generators and rational homotopy
groups therefore gives
\begin{equation}\label{geo:connectivity}
 \pi_j(V_{q,K})\otimes\Q=0\quad(2\leq j<p),
 \qquad
 \pi_p(V_{q,K})\otimes\Q\cong\Q.
\end{equation}
Increasing $K$ makes $p$ arbitrarily large while preserving all the
estimates in \eqref{geo:invariants}.
\end{proof}

The connectivity in the corollary is rational connectivity.  The
manifolds themselves are simply connected, but no vanishing assertion
is made about higher torsion homotopy groups.
Likewise, vanishing rational Pontryagin classes do not imply integral
stable parallelizability.  The construction specifies the complete
rational homotopy type and the indicated characteristic classes;
it does not determine a diffeomorphism type or an explicit
description by bundles or handles.

For a finite-type minimal model with generators in degrees at least
$p\geq2$, cohomology vanishing above $N$, and $\nilh=n\geq1$,
\cite[Theorem~1.1, p.~2]{ParentDC2026} gives the sufficient equality
criterion $N\leq p(n+3)-3\Rightarrow\Hnil=n$.
As a consistency check, our family lies outside this range: with
$n=3q-2$ and $\eta=|a|$,
\begin{equation}\label{geo:dimension-threshold}
 D_{q,K}-\bigl(p(n+3)-3\bigr)
 =(q-3)(q+1)p+B_q(\eta-1)+5-q
 \geq3q^2-4q-1>0.
\end{equation}
The final estimate uses $p\geq2$ and $\eta\geq3$. This comparison
establishes neither a converse to the criterion nor the optimality of
its sufficient threshold.

\section{Spherical stabilization and single-cone attachments}
\label{sec:stabilization}

We now derive a cone-attachment consequence under the spherical
stabilization statement of~\cite{Parent2026}. Coning off the spherical
wedge supplied by that statement recovers the rational homotopy type of
one of our closed, rationally elliptic manifolds. Homology nilpotency
can thereby increase by an arbitrarily large amount.

\Needspace*{12\baselineskip}
The external input is \cite[Theorem~1.2, p.~8; Theorem~15.4 and
Corollary~15.6, pp.~130--131]{Parent2026}, recalled here for convenience.
\nopagebreak
\begin{theorem}[Spherical stabilization, Parent]
\label{stab:quoted}
Let $X$ be a simply connected rational space of finite type, and suppose that
$n=\catz(X)<\infty$. There is a wedge
\[
 S\simeq_{\Q}\bigvee_{\lambda\in L}S_{\Q}^{d_\lambda},
 \qquad d_\lambda\geq2,
\]
with finitely many summands in each dimension, such that
\[
 \catz(X\vee S)=\nilh(X\vee S)=\Hnil(X\vee S)=n.
\]
\end{theorem}

The construction starts from a section of an optimal Sullivan model
for rational category. If $G$ denotes the homotopy
fibre of that section, the associated fibre inclusion $G\to X$ is
nullhomotopic, and its cofibre is rationally $X\vee\Sigma G$. The theorem
shows that $\Sigma G$ is an admissible spherical wedge and that the minimal
model of this cofibre has a quasi-isomorphic quotient of nilpotency $n$.
The next consequence is conditional on this external stabilization
statement; the algebraic and smooth realization theorems of the present
paper do not depend on it.

\begin{theorem}[Unbounded increase under one cone attachment]
\label{stab:cone}
Assume the spherical stabilization statement of Theorem~\ref{stab:quoted}.
Let $q\geq5$ and let $K$ be admissible. There are a simply connected rational
space $Y_{q,K}$ of finite type, a degreewise finite wedge $S_{q,K}$ of
simply connected rational spheres, and a cofibre sequence
\[
 S_{q,K}\xrightarrow{\ i\ }Y_{q,K}
 \longrightarrow (V_{q,K})_{\Q}
\]
such that
\begin{align*}
 Y_{q,K}&\simeq_{\Q}(V_{q,K})_{\Q}\vee S_{q,K},\\
 \Hnil(Y_{q,K})&=\catz(Y_{q,K})=3q-2,\\
 \catz((V_{q,K})_{\Q})&=\catz(Y_{q,K}),\\
 q-4&\leq
 \Hnil((V_{q,K})_{\Q})-\Hnil(Y_{q,K})
 \leq q-2.
\end{align*}
In particular, for every integer $R\geq1$, a single cone attachment can
increase homology nilpotency by at least $R$, with no change in rational
Lusternik--Schnirelmann category.
\end{theorem}

\begin{proof}
Put $X=(V_{q,K})_{\Q}$. Theorems~\ref{alg:main} and~\ref{geo:realization}
give
\[
 \catz(X)=\nilh(X)=3q-2,
 \qquad 4q-6\leq\Hnil(X)\leq4q-4.
\]
Apply Theorem~\ref{stab:quoted} to obtain $S=S_{q,K}$ and set
$Y_{q,K}=X\vee S$. Take $i$ to be the inclusion of the second wedge
summand. Its homotopy cofibre is
\[
 Y_{q,K}\cup_i CS\simeq X,
\]
because coning off the entire summand $S$ leaves the other wedge summand.
The stabilization theorem gives
$\Hnil(Y_{q,K})=\catz(Y_{q,K})=3q-2$.
Subtracting this value from the two bounds for $\Hnil(X)$ proves the
claimed interval. For the last assertion, take $q=R+4$.
\end{proof}

For comparison, the normalized LS category satisfies the classical cone
inequality
\[
 \operatorname{cat}(Z\cup_f CA)\leq\operatorname{cat}(Z)+1,
\]
and the corresponding inequality holds rationally;
see~\cite[Theorem~1.32, p.~16]{CLOT2003}.
Under the stabilization statement, Theorem~\ref{stab:cone} rules out an
analogous uniform additive bound for $\Hnil$. Rational category remains
unchanged in these examples, while the increase in homology nilpotency is unbounded.

\Needspace*{9\baselineskip}
\begin{remark}[The meaning of one cone]
The wedge supplied by Theorem~\ref{stab:quoted} is finite in each dimension;
the theorem does not assert that it has finitely many summands. Accordingly,
Theorem~\ref{stab:cone} attaches one cone on the whole space $S_{q,K}$.
It does not assert that a single cell suffices, or that the source of the
attachment is a finite CW complex. The ellipticity assertion concerns the
manifold $V_{q,K}$ and its rationalization; no ellipticity assertion is made
for the stabilized space $Y_{q,K}$.
\end{remark}

There is also a formulation using ordinary CW complexes. Choose an ordinary
wedge $\widehat S_{q,K}=\bigvee_\lambda S^{d_\lambda}$ with the same sphere
dimensions, and put
$\widehat Y_{q,K}=V_{q,K}\vee\widehat S_{q,K}$. Rationalization gives the
rational wedge in Theorem~\ref{stab:cone}. Thus the spaces in the ordinary
cofibre sequence
\[
 \widehat S_{q,K}\longrightarrow\widehat Y_{q,K}
 \longrightarrow V_{q,K}
\]
have the same values of $\Hnil$ and $\catz$ as their rational counterparts
above. The terminal space is now the closed manifold itself; the invariants
$\Hnil$ and $\catz$ are still rational invariants.

\section{Further consequences and remaining questions}\label{sec:questions}

For the manifolds themselves, the bounds of Theorem~\ref{thm:main}
determine the asymptotic growth of homology nilpotency. For any sequence
of admissible integers $K_q$,
\[
 \lim_{q\to\infty}\frac{\Hnil(V_{q,K_q})}{q}=4,
 \qquad
 \lim_{q\to\infty}
 \frac{\Hnil(V_{q,K_q})}{\catz(V_{q,K_q})}=\frac43.
\]
They leave three possible values for each individual
$\Hnil(V_{q,K})$. Determining which of $4q-6$, $4q-5$,
and $4q-4$ occurs requires either a shorter quasi-isomorphic quotient
or a stronger obstruction than the multiplier used here.

The rational homotopy rank is $q+4$, the category is $3q-2$,
and the least quotient length grows as $4q$. The large degrees
prevent mixing of weights in the obstruction range. Optimizing
the sufficient degree bounds, and hence the realized dimensions,
remains a separate problem.

For the original degrees, what is $\Hnil(M_q)$ as a function of $q$?
In particular, is $\Hnil(M_q)-(3q-2)$ unbounded?
The calculation $\nilh(M_q)=3q-2$ and the upper bound
$\Hnil(M_q)\leq4q-4$ hold without regrading, and
Corollary~\ref{geo:all-regradings} gives smooth realizations in dimension
$d_q$. The lower bound proved here for all graded ideals uses degree
separation and cannot simply be transferred to those original models.

Two further questions concern the size of the gap and the stabilization.
In the present family, $\nilh=3q-2$ grows with $q$. Can the gap be
unbounded while $\nilh$ remains fixed? Can the degreewise finite
spherical stabilization be replaced by a finite wedge, or the cone
attachment by one cell? Resolving these questions requires arguments
beyond the gap estimates and the stabilization result used here.

\appendix

\section{The earlier separation at length three}
\label{app:first-separator}

The first example in~\cite{Parent2026} separating the two invariants
on a fixed minimal model satisfies
\[
 \nilh(M)=3<\Hnil(M)=4.
\]
We recall the construction to place the present family in its historical
context. The separation is proved in
\cite[Theorem~12.4 and Corollary~12.5, p.~100]{Parent2026}, with the
matching upper bound established in
\cite[Theorem~B.2, pp.~149--150]{Parent2026}.

We write $\rho_{16}$ for the generator denoted $r_{16}$
in~\cite{Parent2026}. Let
\[
 B=\frac{\Lambda(s_3,x_5,y_3,a_7,t_{12},\rho_{16})}
 {\bigl(\Lambda^+(s,x,y,a,t,\rho)\bigr)^4}
\]
and prescribe
\[
 ds=dx=dy=0,\qquad da=xy,\qquad dt=ays,
 \qquad d\rho=xt.
\]
These formulas define a differential on the quotient: before truncation,
$d^2\rho=-xays$, which is a word of length four. The cdga $B$ has dimension
$56$ and $(B^+)^4=0$.
Its minimal model can be chosen in the form
\[
 M=\bigl(\Lambda(s_3,x_5,y_3,a_7,t_{12})\otimes\Lambda W,d\bigr)
 \xrightarrow{\simeq}B,
 \qquad W^{\leq18}=0.
\]
The low-degree calculation in~\cite[Corollary~11.8, p.~96]{Parent2026}
then gives
\[
 M^{17}=\Q\{xt\},\qquad Z^{17}(M)=0,
 \qquad d(xt)=-xays\in\m_M^4.
\]
Moreover, $c=axt$ is a cocycle with $[c]\ne0$ in $H^{24}(M)$
\cite[Proposition~11.12, p.~97]{Parent2026}.

These statements isolate the obstruction. If an acyclic differential ideal
$J\subset M$ contained $\m_M^4$, acyclicity would force a primitive of
$d(xt)$ into $J$. The displayed degree calculation makes that primitive
uniquely $xt$. Since $J$ is an ideal and $axt$ is closed,
acyclicity gives $axt\in Z^{24}(J)=d(J^{23})$, contradicting
$[axt]\ne0$ in $H^{24}(M)$. Hence $\Hnil(M)>3$.
On the other hand, $B$ has nilpotency three and $axt$ is an essential
word of length three, so
$e_0(M)=\catz(M)=\nilh(M)=3$. The acyclic ideal constructed in
\cite[Theorem~B.2]{Parent2026} contains $\m_M^5$, proving $\Hnil(M)=4$.

Although $H(M)$ is finite-dimensional, this example is not elliptic.
Indeed, $H^{48}(B)=\Q[\rho^3]$ and
$H^3(B)=\Q[s]\oplus\Q[y]$, whereas $B^{45}=0$.
Thus its cohomology does not satisfy Poincar\'e~duality, which holds for
elliptic minimal Sullivan algebras. This distinguishes the earlier separator
from the finite-generator family constructed in the present paper.

\Needspace*{5\baselineskip}
For this example, \cite[Theorem~15.14, p.~137]{Parent2026}
gives an explicit stabilization by one sphere:
\nopagebreak
\[
 \Hnil\bigl(|M|\vee S_{\Q}^{17}\bigr)=3,
 \qquad\Hnil(|M|)=4.
\]
Coning off this sphere therefore gives an increase of one under a single
rational cell attachment. Under the stabilization statement recalled there,
Section~\ref{sec:stabilization} gives unbounded increases under one cone
attachment, with a possibly infinite, degreewise finite spherical source.

\Needspace*{9\baselineskip}
\begingroup
\makeatletter
\let\addcontentsline\@gobblethree
\section*{Declaration of generative AI and AI-assisted technologies
in the manuscript preparation process}
\makeatother
\endgroup
During the preparation of this manuscript, the author used ChatGPT and
Codex (OpenAI) to assist with literature searches, language editing,
\LaTeX{} preparation, and the development and review of mathematical
arguments. This work was carried out under the author's direction.
The author takes full responsibility for all statements, proofs,
references, and conclusions, and for the article as a whole.

\end{document}